\documentclass{amsart}

\usepackage[latin1]{inputenc}
\usepackage{amssymb,amsmath}
\usepackage{latexsym}
\usepackage{longtable}

\newtheorem{theorem}{Theorem}[section]
\newtheorem{lemma}[theorem]{Lemma}

\newtheorem{corollary}[theorem]{Corollary}
\newtheorem{proposition}[theorem]{Proposition}

\newcommand{\KK}{\overline{K}}
\newcommand{\X}{\overline{X}}
\newcommand{\Y}{\overline{Y}}

\newcommand{\divi}{\operatorname{div}}

\newcommand{\C}{\mathcal C}
\newcommand{\q}{\theta}
\newcommand{\w}{\omega}
\newcommand{\e}{\epsilon}
\newcommand{\HH}{\mathbb H}

\newcommand{\Z}{\mathbb Z}

\newcommand{\CC}{\mathbb C}
\newcommand{\PP}{\mathbb P}

\begin{document}

\title{On the Torelli problem and Jacobian  Nullwerte in genus three }

\author{Jordi Gu\`{a}rdia }

\thanks{Partially supported by  MTM2006-15038-C02-02.}

\address{Departament de Matem\`{a}tica Aplicada IV\\  Escola Polit\`{e}cnica Superior d'Enginyeria de Vilanova i la Geltr\'{u} - Universitat Polit\`{e}cnica de Catalunya\\
Av. V\'{\i}ctor Balaguer s/n. E-08800 Vilanova i la Geltr\'{u}}
\email{guardia@ma4.upc.edu}
\date{\today}


\maketitle

\begin{abstract}

We give a closed formula for recovering a non-hyperelliptic genus three curve from its period matrix, and derive some identities between Jacobian Nullwerte in dimension three.
 \end{abstract}

\section{Introduction}


It is known  that the set of bitangent lines of a non-hyperelliptic genus three curve determines completely the curve, since it admits a unique symplectic structure  (\cite{Caporaso}, \cite{Lehavi}). Given this structure, one can recover an equation for the curve following a method of Riemann (\cite{Riemann},\cite{Ritzenthaler}): one takes an Aronhold system of bitangent lines, determines some parameters by means of some linear systems and then writes down a {\em Riemann model} of the curve. Unfortunately, the parameters involved in this construction
are not defined in general over the field of definition of the curve, but over the field of definition of the bitangent lines.
This is rather inconvenient for arithmetical applications concerned with rationality questions (cf. \cite{Oyono} for instance).

We propose an alternative construction, giving  a model of the curve directly from a certain set of bitangent lines; this model is already defined over the field of definition of the curve. In the particular  case of complex curves, our construction  provides a closed solution for the non-hyperelliptic Torelli problem on genus three:

\begin{theorem} \label{torelli}

Let $\C$ will be a non-hyperelliptic genus three curve defined over a field $K\subset \CC$, and let  $\w_1, \w_2,\w_3$  be a $K$-basis of $ H^0(\C,\Omega^1_{/K}) $, and   $\gamma_1,\dots,\gamma_6$  a symplectic basis of   $H_1(\C,\Z)$.
We denote by $\Omega=(\Omega_1 |\Omega_2)=(\int_{\gamma_j}\omega_k)_{j,k}$ the period matrix of $\C$ with respect to this bases and by
$Z=\Omega_1^{-1}.\Omega_2$ the normalized period matrix.
A model of $\C$ defined (up to normalization) over $K$ is:
$$
\begin{array}{c}\displaystyle
\left(
\frac{[w_7w_2w_3][w_7w'_2w'_3]}{[w_1w_2w_3][w'_1w'_2w'_3]}
X_1Y_1
  +
\frac{[w_1w_7w_3][w'_1w_7w'_3]}{[w_1w_2w_3][w'_1w'_2w'_3]}X_2Y_2   -
\frac{[w_1w_2w_7][w'_1w'_2w_7]}{[w_1w_2w_3][w'_1w'_2w'_3]} X_3Y_3
\right)^2\\\\
\displaystyle
-4
\frac{[w_7w_2w_3][w_7w'_2w'_3]}{[w_1w_2w_3][w'_1w'_2w'_3]}\frac{[w_1w_7w_3][w'_1w_7w'_3]}{[w_1w_2w_3][w'_1w'_2w'_3]}
X_1 Y_1 X_2 Y_2=0,
\end{array}
$$
where
$$
\begin{array}{ll}
w_1={}^t(\frac12,0,\frac12)+{}^t(0,0,\frac12).Z &
w_1'={}^t(0,0,\frac12)+{}^t(0,0,\frac12).Z \\\\
w_2={}^t(\frac12,\frac12,0)+{}^t(0,\frac12,\frac12).Z &
w_2'={}^t(0,\frac12,0)+{}^t(0,\frac12,\frac12).Z \\\\
w_3={}^t(\frac12,\frac12,\frac12)+{}^t(0,\frac12,0).Z &
w_3'={}^t(0,\frac12,\frac12)+{}^t(0,\frac12,0).Z \\\\
w_7={}^t(0,0,\frac12)+{}^t(\frac12,\frac12,\frac12).Z.
\\\\
\hskip 1cm X_j=\operatorname{grad}_{z=0} \q[w_j](z;Z) .\Omega_1^{-1}. (X,Y,Z)^{t},\\\\
\hskip 1cm  Y_j=\operatorname{grad}_{z=0} \q[w'_j](z;Z) .\Omega_1^{-1}.(X,Y,Z)^{t},
\end{array}
$$
and $[u,v,w]$ denotes the Jacobian Nullwert given by $u, v, w$.
\end{theorem}

Of course, one could simplify the denominators in the equation above to obtain a simpler formula. The advantage of writing the fractions is  that their values are algebraic over the field of definition of the curve.

The formula above can be interpreted as a universal curve over the moduli space of complex non-hyperelliptic curves of genus three. Moreover, using  the Frobenius formula, we can express the equation of the curve in terms of Thetanullwerte instead of Jacobian Nullwerte (cf. section \ref{frobenius}).

As it will become clear along the paper, many different choices for the $w_k$ are possible.
They are only restricted to some geometric conditions expressed in terms of the associated bitangent lines (cf.  section \ref{teoremes}). Every election of the $w_k$ leads to a model for the curve, though all these models agree up to a proportionality constant.
  The comparison of the models given for different elections provides a number of identities between Jacobian Nullwerte in dimension three. For instance, we prove:

\begin{theorem}\label{igualtats} Take $w_1, w_2, w_3, w'_1, w_2', w_3'$ as before, and write
$$
\begin{array}{ll}
w_4={}^t(\frac12,0,0)+{}^t(\frac12,\frac12,\frac12).Z, 
&
w_4'={}^t(0,0,0)+{}^t(\frac12,\frac12,\frac12),\\\\
w_7'={}^t(\frac12,\frac12,0)+{}^t(\frac12,0,\frac12).Z.
\end{array}
$$
 The following equalities hold on the space $\HH_3$:
$$
\begin{array}{c}
\displaystyle
{[ w_2 w_3 w_7](Z)}{ [w_1 w_3' w'_7 ](Z)}=    {[ w_1 w_3 w_7](Z)} {[w_2 w_3' w'_7 ](Z)},
\\\\
{[ w'_2 w_3 w'_7](Z)}{ [w'_1 w'_3 w_7 ](Z)}  ={[w'_1 w_3 w_7'](Z)}  {[w_2' w'_3 w_7 ](Z)},\\\\ \\\\
\displaystyle
{[w_3 w_1 w_2](Z)[w'_3 w_1 w_2](Z)}
{[w_4 w_1 w'_2](Z)[w'_4 w_1 w'_2](Z) }\\ \shortparallel\\
 {[w_4 w_1 w_2](Z)[w'_4 w_1 w_2](Z)}{[w_3 w_1 w'_2](Z)[w'_3 w_1 w'_2](Z) },\\\\\\
\end{array}
$$
$$
\begin{array}{c}
\displaystyle
{[w_3 w_1 w_2](Z)[w'_3 w_1 w_2](Z)}{[w_4 w_2 w'_1](Z)[w'_4 w_2 w'_1](Z) }\\\shortparallel\\
 {[w_4 w_1 w_2](Z)[w'_4 w_1 w_2](Z)}
{[w_3 w_2 w'_1](Z)[w'_3 w_2 w'_1](Z) },\\\\\\
\displaystyle
{[w_3 w_1 w_2](Z)[w'_3 w_1 w_2](Z)}{[w_4 w'_1 w'_2](Z)[w'_4 w'_1 w'_2](Z) }\\\shortparallel\\
{[w_4 w_1 w_2](Z)[w'_4 w_1 w_2](Z)}{[w_3 w'_1 w'_2](Z)[w'_3 w'_1 w'_2](Z) }.
\end{array}
$$
\end{theorem}

Again, there are many possible elections of the $w_k$, each one leading to  a set of  identities between Jacobian Nullwerte.

Apart from its intrinsical theoretical interest, these results have different applications. From the computational viewpoint, for instance, theorem \ref{torelli} (or corollary \ref{torelli-theta}) can be used to determine equations for modular curves or to present three-dimensional factors of modular jacobians as jacobians of curves (thus improving the results in \cite{Oyono}). In a more theoretical frame, the identities in (\ref{igualtats}) may lead to simplified expressions for the discriminant of genus three curves (cf. {\cite{Ritzenthaler}).

The results stated above are the analytic version of the corresponding results in an algebraic context. These results are proved in the first part of the paper. We start recalling Riemann construction for curves of genus three, and the basic concepts regarding the symplectic structure of the set of bitangent lines of these curves. The relation between different Steiner complexes is described explicitly in section 3. In sections 4 and 5 we combine the ideas in previous sections to provide the algebraic versions of our theorems. The second part of the paper contains the results in the analytic context. The proof of theorems \ref{teo-formula} and \ref{igualtats} is given in section 6. Finally, we give the Thetanullwerte version of theorem \ref{teo-formula} in section 7, which also contains a geometric description of the fundamental systems appearing in the Frobenius formula.

The author is much indebted to C. Ritzenthaler for pointing his attention to the Torelli problem for non-hyperelliptic genus three curves, and to E. Nart for many fruitful discussions.

$\,$
%

{\bf Notation and conventions:} We will work with a  non-singular genus three curve $\C$, defined over a field $K$ with $\operatorname{char}(K)\ne$2, and given by a quartic  equation $Q=0$. We assume that the curve is embedded in the projective plane $\mathbb{P}^2(K)$ by its canonical map. Given two polynomials $Q_1, Q_2$, we will write $Q_1\sim Q_2$ to express that they agree up to a constant proportionality factor.
 A number of geometric concepts and facts about genus three curves will be used along the paper. We will take \cite{Dolgachev}  as basic reference, following the notations introduced there.


\part{ALGEBRAIC IDENTITIES}

\section{Riemann model and Steiner complexes}

A triplet or tetrad of bitangent lines is called syzygetic whenever their contact points with the curve lie on a conic.
Any pair of bitangent lines to $\C$ can be completed in five  different ways with another pair to a syzygetic tetrad. Moreover, any two of the six pairs form a syzygetic tetrad. Such a set of  of six pairs of bitangent lines  is called {\em Steiner complex} (cf. \cite{Dolgachev}).

Riemann showed that, given three pairs of bitangent lines in a Steiner complex, one can determine proper equations  $\{X_1=0,Y_1=0\},\{X_2=0,Y_2=0\}$,  $\{X_3=0,Y_3=0\}$   of these lines so that
the quartic equation  $Q=0$ defining  the  curve $\C$ can be written as:
$$
Q\sim( X_1 Y_1+ X_2 Y_2-X_3Y_3)^2-4 X_1Y_1 X_2Y_2=0,
$$
 and moreover they satisfy the relation:
$$
X_1+X_2+X_3=Y_1+Y_2+Y_3. 
$$
To find the proper equations of the bitangent lines, one takes arbitrary equations for them and solves certain linear systems to compute  scaling factors leading to the well adjusted  equations.

Reciprocally, whenever we take three pairs  $\{X_1, Y_1\},\{X_2, Y_2\} \{X_3, Y_3\}$  of linear polynomials over $K$, such that no triplet formed by a polynomial in each pair is linearly dependent and $X_1+X_2+X_3=Y_1+Y_2+Y_3$, the equation above gives a non-singular genus three curve over $K$, with   $\{X_1=0, Y_1=0\}$, $\{X_2=0, Y_2=0\}$, $\{X_3=0, Y_3=0\}$ being three pairs of bitangent lines in a common Steiner complex.

Let us define $W=X_1+X_2+X_3=0$, $Z_1=X_1+Y_2+Y_3$, $Z_2=Y_1+X_2+Y_3$, $Z_3=Y_1+Y_2+X_3$. The trivial equalities:
\begin{equation}
\begin{array}{rl}
Q
\sim&\left(X_1 Y_2+X_2 Y_1-W Z_3\right)^2-4 X_1 X_2 Y_1 Y_2\\
=&(X_1 Y_3+X_3 Y_1- W Z_2)^2-4 X_1 X_3 Y_1 Y_3\\
=&(X_2 Y_3+X_3 Y_2-W Z_1)^2-4X_2 X_3 Y_2 Y_3 \\\\
=& (X_1 X_2 + Y_1 Y_2-Z_1 Z_2)^2-4 X_1 X_2 Y_1 Y_2\\
=& (X_1 X_3+Y_1 Y_3-Z_1 Z_3)^2-4 X_1 X_3 Y_1 Y_3\\
=& (X_2 X_3+Y_2 Y_3-Z_2 Z_3)^2-4 X_1 X_3 Y_1 Y_3\\\\
=&(X_1 W+Y_2 Z_3-Y_3 Z_2)^2-4X_1 W Y_2 Z_3\\
=&(X_2 W+Y_1 Z_3-Y_3 Z_1)^2-4X_2 W Y_1 Z_3\\
=&(X_3 W+Y_1 Z_2-Y_2 Z_1)^2-4X_3 W Y_1 Z_2\\\\
\end{array}
$$
$$
\begin{array}{rl}
=&(X_1 Z_3+Y_2 W-X_3 Z_1)^2-4X_1 Z_3 Y_2 W  \\
=&(X_1 Z_2+Y_3 W-X_2 Z_1)^2-4X_1 Z_2 Y_3 W  \\\\
=&(X_1 Z_1+X_2 Z_2-X_3 Z_3)^2-4X_1 Z_1 X_2 Z_2  \\
=&(Y_1 Z_1+Y_2 Z_2-Y_3 Z_3)^2-4X_1 Y_1 Y_2 Z_2      \\\\
\end{array}
\end{equation}
show that $W=0, Z_1=0, Z_2=0, Z_3=0$ are also bitangent lines to $\C$, and they
make apparent a number of different Steiner complexes on $\C$.

On the other hand, one can check easily that in this situation
any triplet formed picking a line from each pair $\{X_i=0, Y_i=0\}$ is asyzygetic. Analogous conclusions can be derived from each of the equations above. For instance, we mention for our later convenience that   $\{X_1=0, X_2=0, W=0\}$, $\{X_1=0, X_3=0, W=0\}$, $\{ X_1=0, Y_2=0, Z_3=0\}$ and $\{X_1=0, Y_3=0, Z_3=0\}$ are asyzygetic triplets.

\section{Relations between Steiner complexes}

Every Steiner complex $S$ has an associated two-torsion element  $O(D_S)\in\operatorname{Pic}^0(\C)$; a divisor $D_S$ defining it is given by the difference of the contact points of the two bitangent lines in any pair in $S$. We will call such a divisor an {\em associated divisor} of $S$.
Indeed, the map $S\mapsto O(D_S)$ establishes a bijection between the set of Steiner complexes on $\C$ and $\operatorname{Pic}^0(\C)[2]$  (cf. \cite{Dolgachev}).

Two different Steiner complexes share four or six bitangent lines; they are called syzygetic or asyzygetic accordingly. There is a simple criterion to check whether two Steiner complexes are syzygetic:

\begin{lemma}{\cite[p. 96]{Dolgachev}}\label{dolgachev}
Let $S_1, S_2$ be two Steiner complexes, and let  $D_1, D_2$ be associated divisors to them. Then $\sharp \overline{S_1}\cap \overline{S_2}=4$ if
$e_2(D_1, D_2)=0$ and $\sharp \overline{S_1}\cap \overline{S_2}=6$ if $e_2(D_1, D_2)=1$, where $e_2$ denotes the Weil pairing on $\operatorname{Pic}^0(\C)[2]$.
\end{lemma}

We note that the Weil pairing can be computed by means of the
Riemann-Mumford relation (\cite[p.290]{ACGH84}): for any semicanonical divisor $D$, we have
\begin{equation}\label{Riemann-Mumford}
e_2(D_1,D_2)=h^0(D)+h^0(D+D_1+D_2)+h^0(D+D_1)+h^0(D+D_2)   \pmod2.
\end{equation}

Given a pair   $\{X=0,
Y=0\}$ of bitangent lines, we shall denote by
 $S_{XY}$  the Steiner complex determined by them. For any set $S=\{\{X_1=0, Y_1=0\}, \dots, \{X_r=0, Y_r=0\}\}$ of pairs of bitangent lines, the subjacent set of lines
will be denoted by $\overline{S}=\{X_1=0, Y_1=0, \dots, X_r=0, Y_r=0\}$.

A priori,
the bitangent lines which form a pair in a Steiner complex, are completely indistinguishable. But if we consider the Steiner complex with relation to others, there appears a individualization of every line. This idea is made explicit in  corollary \ref{tripletes-syzyg} and proposition \ref{steiner-asyzygetics}.

It is evident that whenever  $ \{X_1=0, Y_1=0\}, \{X_2=0, Y_2=0\}$
is a Steiner couple,   the pairs  $\{X_1=0, Y_2=0\}, \{X_2=0,
Y_1=0\}$ form also a Steiner couple, and also  the pairs
$\{X_1=0, X_2=0\}, \{Y_1=0,
Y_2=0\}$ do.
The corresponding Steiner complexes are tightly related:

\begin{proposition}
Let $S_{X_1Y_1}=\{ \{X_1=0, Y_1=0\},\dots, \{X_6=0, Y_6=0\}\}$ be a
Steiner complex.
\begin{itemize}
\item[a)] The Steiner complexes  $S_{X_1Y_1}$ and $S_{X_1Y_j}$ are syzygetic and they share the four lines  $X_1=0, Y_1=0, X_j=0,
Y_j=0$, i.e., $\overline{S}_{X_1 Y_1}\cap \overline{S}_{X_1
Y_j}=\left\{X_1=0,\right.$ $\left. Y_1=0, X_j=0, Y_j=0\right\}$.
\item[b)] The Steiner complexes $S_{X_1Y_1}$ and $S_{X_1 X_2}$ are syzygetic and $\overline{S}_{X_1 Y_1}\cap \overline{S}_{X_1
X_2}=\{X_1=0,\, Y_1=0, X_2=0, Y_2=0\}$.
\item[c)] The Steiner complexes $S_{X_1X_2}$ and $S_{X_1 Y_2}$ are syzygetic and $\overline{S}_{X_1 Y_1}\cap \overline{S}_{X_1
X_2}=\{X_1=0,\, Y_1=0, X_2=0, Y_2=0\}$.
\item[d)] For $j\ne k$, $j,k\ne 1$ the Steiner complexes $S_{X_1Y_j}$ and $S_{X_1Y_k}$ are
asyzygetic.
\end{itemize}
\end{proposition}

\noindent {\bf Proof:} Let us write $\divi(X_i):=2P_i+2Q_i$, $\divi(Y_i):=2R_i+2S_i$. We apply the criterion of lemma \ref{dolgachev}, using formula \ref{Riemann-Mumford} to
compute the Weil pairing of  divisors $D_{11}$ and $D_{1j}$
associated to $S_{X_1Y_1}$ and $S_{X_1Y_j}$ respectively.
We take $D=P_1+Q_1$ and find
$$
\begin{array}{rl}
e_2(D_{11},D_{1j})&=h^0(3P_1+3Q_1-R_1-S_1-R_j-S_j)-h^0(2P_1+2Q_1-R_1-S_1)\\
&-h^0(2P_1+2Q_1-R_j-S_j)+h^0(P_1+Q_1)=\\
=&h^0(K_{\C}+P_1+Q_1-R_1-S_1-R_j-S_j)-h^0(K_{\C}-R_1-S_1)\\
&-h^0(K_{\C}-R_j-S_j)+1=\\
=&h^0(P_1+Q_1+R_1+S_1-R_j-S_j)-h^0(R_1+S_1)-h^0(R_j+S_j)+1\\
=& h^0(P_1+Q_1+R_1+S_1-R_j-S_j)-1=0,
\end{array}
$$
since $\{X_1=0, Y_1=0, X_j=0,Y_j=0\}$ form a syzygetic tetrad. This proves part {\it a}). Remaining parts  are proved analogously.
$\blacksquare$\\

We derive from here a more explicit version of theorem 6.1.8 in \cite{Dolgachev}:

\begin{corollary} \label{total}
Let $\{ \{X_1=0, Y_1=0\}, \{X_2=0, Y_2=0\}\}$ be a syzygetic tetrad of bitangent lines..
The Steiner complexes $S_{X_1 Y_1}$, $S_{X_1 Y_2}$, $S_{X_1 X_2}$
satisfy:
$$
\overline{S}_{X_1Y_1}\cup
\overline{S}_{X_1Y_2}\cup\overline{S}_{X_1X_2}=\operatorname{Bit}(\C).
$$
\end{corollary}

Some immediate consequences of the above proposition are:

\begin{corollary} \label{tripletes-syzyg}
Let $S=\{ \{X_1=0, Y_1=0\},\dots, \{X_6=0, Y_6=0\}\}$ be a Steiner complex.
\begin{itemize}
\item[a)]
Any triplet   $\{X_i=0, Y_i=0, X_j=0\}$  is syzygetic.
\item[b)]
Any triplet   $\{X_i=0, X_j=0, X_k=0 \}$ formed picking a line  from three different pairs of  $S$ is asyzygetic.
\item[c)] Let $\{U=0, V = 0 \}$ be an arbitrary pair in the Steiner complex $S_{X_1 Y_2}$. Then
$U=0\in$ belongs to $S_{X_1 X_3}$ but not to $S_{X_1 Y_3}$ and $V=0$ belongs to $S_{X_1 Y_3}$ but not to $S_{X_1 X_3}$ (or the other way round).
\end{itemize}
\end{corollary}

We now derive a geometrical property of asyzygetic triplets of bitangent lines that we will need later.

\begin{corollary}
Every three asyzygetic bitangent lines $X_1=0, X_2=0, X_3=0$ can be paired with other three bitangent lines $Y_1=0$, $Y_2=0$, $Y_3=0$ so that the three pairs $\{X_i=0, Y_i=0\}$ belong to a common Steiner complex.
\end{corollary}


\begin{proposition}\label{noconcorrents} Three asyzygetic bitangent  lines do not cross in a point.
\end{proposition}

\noindent {\bf Proof:} Complete the lines to three pairs  $\{ \{X_1=0, Y_1=0\}, \{X_2=0, Y_2=0\}, \{X_3=0, Y_3=0\}\}$ in a common Steiner complex, and re-scale the equations to have a Riemann model for $\C$:
$
(X_1 Y_1+X_2 Y_2-X_3Y_3)^2-4 X_1
Y_1 X_2Y_2=0.
$
If the lines $X_1=0$, $X_2=0$, $X_3=0$ crossed in a point, this
should be a singular point of $\C$, which is non-singular by
hypothesis. $\blacksquare$\\

We end this section with a explicit description of  triplets of mutually asyzygetic Steiner complexes:

\begin{proposition} \label{steiner-asyzygetics}
Let $X_1=0, X_2=0, X_3=0$ be three asyzygetic bitangent lines. After
a proper labeling of the bitangent lines, the Steiner complexes
$S_{X_1 X_2}$, $S_{X_2,X_3}$ and $S_{X_3 X_1}$ have the following
shape:
$$
S_{X_1 X_2}=\begin{pmatrix} X_1 & X_2 \\ X_4 & X_9 \\ X_5 & X_{10}\\
X_6 & X_{11} \\ X_7 & X_{12}\\ X_8 & X_{13}
\end{pmatrix},\quad
S_{X_2 X_3}=\begin{pmatrix} X_2 & X_3 \\ X_9 & X_{14} \\ X_{10} & X_{15}\\
X_{11} & X_{16} \\ X_{12} & X_{17}\\ X_{13} & X_{18}
\end{pmatrix},\quad
S_{X_3 X_1}=\begin{pmatrix} X_3 & X_1 \\  X_{14} & X_4 \\ X_{15} & X_{5}\\
X_{16} & X_{6} \\ X_{17} & X_{7}\\ X_{18} & X_{8}
\end{pmatrix}.
$$
In particular, the three complexes are asyzygetic.
\end{proposition}

\noindent\textbf{Proof:} Take a second pair of bitangent lines $\{A=0, B=0\}$ in  $S_{X_1X_2}$.  After corollary \ref{total}, any other bitangent line must lie in exactly one of the Steiner complexes
$$
\begin{array}{l}
S_{X_1X_2}=\{\{X_1=0, X_2=0\}, \{A=0, B=0\},\dots\},\\
S_{X_1A}=\{\{X_1=0, A=0\}, \{X_2=0, B=0\},\dots\},\\
S_{X_1B}=\{\{X_1=0, B=0\}, \{X_2=0, A=0\},\dots\}.
\end{array}
$$
Since $X_3=0$ cannot be in $S_{X_1X_2}$ by hypothesis, we may suppose that we have a pair $\{X_3=0, C=0\}$ in
$S_{X_1 A}$, so that $X_1=0, A=0, X_3=0, C=0$ is a syzygetic tetrad, and thus $\{A=0, C=0\}$ belongs to $S_{X_1X_3}$; but also $X_2=0, B=0, X_3=0, C=0$  is a syzygetic tetrad and hence $\{B=0, C=0\}$ belongs to $S_{X_2X_3}$.
$\blacksquare$

\section{From bitangent lines to equations}\label{teoremes}

The basic tool for the proof of theorem \ref{torelli} is an algebraic version of it, giving a general equation for a curve of genus three in terms of some of its bitangent lines:

\begin{theorem}\label{teo-formula} Let $\C: Q=0$ be a non-singular genus three plane curve defined over a field $K$ with $\operatorname{char}(K)\ne$2.
Let $\{X_1=0, Y_1=0\}$, $\{X_2=0, Y_2=0\}$, $\{X_3=0, Y_3=0\}$ be
three pairs of bitangent lines from a given Steiner complex. Let
$\{X_7=0, Y_7=0\}$ be a fourth pair of bitangent lines from the Steiner complex given by the pairs $\{\{X_1=0, Y_2=0\}, \{X_2=0,
Y_1=0\}\}$, ordered so that $\{X_1=0, X_3=0, X_7=0\}$ and $\{X_1=0,
Y_3=0, Y_7=0\}$ are asyzygetic.

The equation of $\C$ can be presented as
\begin{equation}\label{formula-det}
\begin{array}{c}
Q\sim \left(\frac{(X_7X_2X_3)(X_7Y_2Y_3)}{(X_1X_2X_3)(Y_1Y_2Y_3)}X_1Y_1+
\frac{(X_1X_7X_3)(Y_1X_7Y_3)}{(X_1X_2X_3)(Y_1Y_2Y_3)}X_2Y_2-
\frac{(X_1X_2X_7)(Y_1Y_2X_7)}{(X_1X_2X_3)(Y_1Y_2Y_3)}X_3Y_3\right)^2\\\\
-4\frac{(X_7X_2X_3)(X_7Y_2Y_3)}{(X_1X_2X_3)(Y_1Y_2Y_3)}\frac{(X_1X_7X_3)(Y_1X_7Y_3)}{(X_1X_2X_3)(Y_1Y_2Y_3)}
X_1 Y_1 X_2 Y_2=0.
\end{array}
\end{equation}

Here $(ABC)$ denotes the determinant of the matrix formed by the
coefficients of the homogeneous linear polynomials $A, B, C$.
\end{theorem}

\noindent
\textbf{Remark:}
All the determinants appearing in the theorem are
non-zero, since they are formed with asyzygetic triplets of
bitangent lines, which are non-concurrent by proposition
\ref{noconcorrents}.

\noindent {\bf Proof:} We know that for a proper re-scaling $\overline{X_i}=\alpha_i X_i$,
$\overline{Y_i}=\beta_i X_i$, we will have the following  Riemann model for $\C$:
$$
Q= ( \X_1 Y_1+ \X_2Y_2- \X_3Y_3)^2-4\X_1Y_1\X_2Y_2=0,
$$
with  $\X_7=\X_1+\X_2+\X_3=-\Y_1-\Y_2-\Y_3=0$. We look at these two equalities as equations in the scaling factors $\alpha_i, \beta_i$ to determine them. We find:
$$
\begin{array}{lll}
\displaystyle \alpha_1=\alpha_7\frac{( X_7X_2X_3)}{(X_1X_2X_3)},\qquad
&
\displaystyle \alpha_2=\alpha_7\frac{(X_1  X_7X_3)}{(X_1X_2X_3)},\qquad
&
\displaystyle \alpha_3=\alpha_7\frac{(X_1X_2 X_7)}{(X_1X_2X_3)}, \\\\
\displaystyle \beta_1 =-\alpha_7\frac{( X_7Y_2Y_3)}{(Y_1Y_2Y_3)},\qquad
&
\displaystyle \beta_2=-\alpha_7\frac{(Y_1 X_7Y_3)}{(Y_1Y_2Y_3)},\qquad
&
\displaystyle \beta_3=-\alpha_7\frac{(Y_1Y_2 X_7)}{(Y_1Y_2Y_3)}.
\end{array}
$$
We substitute the $\X_i, \Y_i$ on the Riemann model by his values, and simplify the constant $\alpha_7$ to get the desired equation. $\blacksquare$


\section{Determinants of bitangent lines}\label{dets}

 We have seen in theorem \ref{teo-formula} how to construct a presentation of the curve $\C$ from certain set of bitangent lines, which we can choose in many different ways. If  we make different  elections, the comparison of the corresponding presentations, which all agree up to a constant, leads us to a number of equalities between  the involved determinants. It turns out that the identities obtained in this way can be easily proved independently of their deduction.

\subsection{Relations between pairs  in the same Steiner complex}

\begin{theorem}\label{dobles}
Let $\{X_1=0, Y_1=0\},\dots, \{X_4=0, Y_4=0\}$ be four different pairs in the Steiner complex $S=S_{X_iY_i}$. Then
$$
\begin{array}{rl}
\displaystyle\frac{(X_3 X_1 X_2)(Y_3 X_1 X_2)}{(X_4 X_1 X_2)(Y_4 X_1 X_2)}
&=\displaystyle\frac{(X_3 X_1 Y_2)(Y_3 X_1 Y_2) }{(X_4 X_1 Y_2)(Y_4 X_1 Y_2) }\\\\
&=\displaystyle\frac{(X_3 X_2 Y_1)(Y_3 X_2 Y_1) }{(X_4 X_2 Y_1)(Y_4 X_2 Y_1) }=\frac{(X_3 Y_1 Y_2)(Y_3 Y_1 Y_2) }{(X_4 Y_1 Y_2)(Y_4 Y_1 Y_2) }.
\end{array}
$$
\end{theorem}

\noindent{\bf Proof:}  It is known (cf. \cite[p. 94]{Dolgachev}) that all the pairs in a Steiner complex
can be seen as one of the  degenerate conics on a fixed pencil of conics. In particular,  we must have a relation
$$
\lambda^2 X_1 Y_1 +X_2 Y_2+ \lambda X_3 Y_3=X_4 Y_4,
$$
for certain $\lambda \in \KK^\ast$.
Let $P_1, P_2, P_3, P_4$ be the points of intersection of the
pairs of lines $\{X_1=X_2=0\}$,  $\{X_1=Y_2=0\}$,  $\{Y_1=X_2=0\}$,
$\{Y_1=Y_2=0\}$.
If we substitute these points in the equality above we find:
$$
\lambda X_3(P_j)Y_3(P_j)=X_4(P_j)Y_4(P_j), j=1,\dots,4,
$$
so that for $j\ne k$ we have
\begin{equation}\label{34}
\frac{X_3(P_j)Y_3(P_j)}{X_4(P_j)Y_4(P_j)}=\frac{X_3(P_k)Y_3(P_k)}{X_4(P_k)Y_4(P_k)}.
\end{equation}
An elementary exercise in linear algebra shows that for any two lines $U=0, V=0$ we have the  equalities
$$
\frac{U(P_1)}{V(P_1)}=\frac{(U X_1 X_2)}{(V X_1 X_2)},
\frac{U(P_2)}{V(P_2)}=\frac{(U X_1 Y_2)}{(V X_1 Y_2)},
\frac{U(P_3)}{V(P_3)}=\frac{(U Y_1 X_2)}{(V Y_1 X_2)},
\frac{U(P_4)}{V(P_4)}=\frac{(U Y_1 Y_2)}{(V Y_1 Y_2)}.
$$
The theorem follows if we apply these relations to the two pairs of lines $X_3=0, X_4=0$ and $\nolinebreak{Y_3=0}$, $\nolinebreak{Y_4=0}$  and substitute the resulting relations in (\ref{34}). $\blacksquare$

\subsection{Relations between pairs  in syzygetic Steiner complexes}

\begin{theorem}\label{simples}
Let $\{X_1=0, Y_1=0\}$, $\{X_2=0, Y_2=0\}$, $\{X_3=0, Y_3=0\}$ be
three pairs of bitangent lines from a given Steiner complex. Let
$\{X_7=0, Y_7=0\}$ be a fourth pair of bitangent lines from the Steiner complex given by the pairs $\{\{X_1=0, Y_2=0\}, \{X_2=0,
Y_1=0\}\}$, ordered so that $\{X_1=0, X_3=0, X_7=0\}$ and $\{X_1=0,
Y_3=0, Y_7=0\}$ are asyzygetic.
The following relations hold:
$$
\frac{(X_2 X_3 X_7)}{(X_1 X_3 X_7)}=\frac{(X_2 Y_3  Y_7)}{(X_1 Y_3 Y_7)},
\qquad
 \frac{(Y_2 X_3  Y_7)}{(Y_1 X_3 Y_7)}=
\frac{(Y_2 Y_3 X_7)}{(Y_1 Y_3 X_7)}.
$$
\end{theorem}

\noindent{\bf Proof:} The validity of the equalities is not affected by a re-scaling of the involved lines, so that we can assume that $X_7=X_1+X_2+X_3$, $Y_7=Y_1+Y_2+X_3=-X_1-X_2-Y_3$.
The result follows from the very elementary properties of the determinants. $\blacksquare$\\

\part{ANALYTIC IDENTITIES}

We  now prove the theorems stated in the introduction, which turn out to be the translation of the results in the first part of the paper to the context of complex geometry.


From now on, we suppose that $\C$ is a non-hyperelliptic plane curve of genus three defined over a field $K\subset \CC$. We choose the basis $\w_1, \w_2,\w_3\in H^0(\C,\Omega^1_{/K}) $ of holomorphic differential forms such that the we can identify $\C$  with the image of the associated canonical map $\iota: \C\longrightarrow \PP^2(K)$. We also fix
 a symplectic basis $\gamma_1,\dots,\gamma_6$ of the singular homology $H_1(\C,\Z)$. We denote by $\Omega=(\Omega_1 |\Omega_2)=(\int_{\gamma_j}\omega_k)_{j,k}$ the period matrix of $\C$ with respect to these bases, and by
$Z=\Omega_1^{-1}.\Omega_2$ the normalized period matrix. We consider the Jacobian variety $J_\C$, represented by the complex torus
$\mathbb{C}^3/(1|Z)$, with the Abel-Jacobi map:
\begin{equation}\label{Abel-Jacobi}
\begin{array}{rcl}
\C^2 &\stackrel{\Pi}\longrightarrow & J_{\C} \\
D & \longrightarrow &\int_\kappa^D(\w_1,\w_2,\w_3),
\end{array}
\end{equation}
where $\kappa$ is the Riemann vector, which guarantees  that $\Theta=\Pi(\C^2)$ is the divisor  cut
by the Riemann theta function $\q(z;Z)$, and that
$\Pi(K_{\C}-D)=-\Pi(D)$.

We shall describe, as customary, the elements of $J_\C[2]$ by means of {\em characteristics}: every $w\in J_\C[2]$ is determined uniquely  by a six-dimensional vector $\e=(\e',\e'')\in \{0,1/2\}^6$ by the relation $w=\e+Z\e''$.
With this notation, the Weil pairing on $J_\C[2]$ is given by:
\begin{equation}\label{Weil}
\tilde e_2(w_1,w_2):=\tilde e_2(\e_1,\e_2):=\e_1'\cdot \e_2''+\e_1''\e_2'\pmod2.
\end{equation}

Every  $w\in J_\C[2] $
defines a translate of the Riemann theta function:
$$
\begin{array}{rl}
\q[w](z;Z):=& e^{\pi i ^t\e'.Z.\e'+2\pi i ^t\e'.(z+\e'')}\q(z+Z\e'+\e'') \\ \\
=&\displaystyle\sum_{n\in \mathbb{Z}^{g}}e^{\pi i^{t}(n+\e').Z.(n+\e')+2\pi i^{t}(n+\e').(z+\e'')}.
\end{array}
$$
The values $\q[w](0,Z)$ are usually called {\em Thetanullwerte}, and
denoted shortly by $\q[w](Z)$. For a sequence of three points $w_1, w_2, w_3\in J_\C[2]$ the modified determinant $[w_1, w_2, w_3](Z):=\pi^3 \det J[w_1,w_2,w_3](Z)$ of the matrix
$$
J[w_1,w_2,w_3](Z):=\left(\begin{array}{ccc}
\displaystyle\frac{\partial \theta [w _{1}]}{\partial z_{1}}(0;Z) &
\displaystyle\frac{\partial \theta [w _{1}]}{\partial z_{2}}(0;Z) &
\displaystyle\frac{\partial \theta [w _{1}]}{\partial z_{3}}(0;Z) \\
\displaystyle\frac{\partial \theta [w _{2}]}{\partial z_{1}}(0;Z) &
\displaystyle\frac{\partial \theta [w _{2}]}{\partial z_{2}}(0;Z) &
\displaystyle\frac{\partial \theta [w _{2}]}{\partial z_{3}}(0;Z)
\\
\displaystyle\frac{\partial \theta [w _{3}]}{\partial z_{1}}(0;Z) &
\displaystyle\frac{\partial \theta [w _{3}]}{\partial z_{2}}(0;Z) &
\displaystyle\frac{\partial \theta [w _{3}]}{\partial z_{3}}(0;Z)
\end{array}
\right)
$$
is called {\em Jacobian Nullwert}.
(These definitions can be given for a $g$-dimensional complex torus abelian varieties, but we restrict them to our case of dimension three for brevity).

\section{Proof of main results}

It is well-known that the Abel-Jacobi map establishes a bijection between semicanonical divisors on $\C$ and the set  $J_\C[2]$.  In our particular case of a non-hyperelliptic genus three curve, the odd semicanonical divisors are given by the bitangent lines: given a bitangent $X=0$, then $\divi(X)=2D$ with $D$ a semicanonical divisor which goes to a $w:=\Pi(D)\in J_\C[2]^{\rm odd}$. Reciprocally, given $w\in J_\C[2]^{\rm odd}$, the line
$$
H_w:= \left(\frac{\partial \theta [w]}{\partial z_{1}}(0;Z),\frac{\partial \theta [w]}{\partial z_{2}}(0;Z),
\frac{\partial \theta [w]}{\partial z_{3}}(0;Z)\right).\Omega_1^{-1}.\left(
\begin{array}{c}
X \\ Y \\ Z
\end{array}
\right)=0,
$$
is a bitangent line to $\C$, whose corresponding semicanonical divisor $D$  satisfies $\Pi(D)=w$.
Moreover, we have
$$
e_2(\Pi(w_1),\Pi(w_2))=\tilde e_2(w_1,w_2).
$$
Note also that, given  $w_i \in J_\C[2]^{\rm odd}$, we have
\begin{equation}\label{trivial}
\displaystyle\frac{[w_1,w_2,w_3](Z)}{[w_4,w_5,w_6](Z)}=\frac{(H_{w_1}  H_{w_2}
H_{w_3})}{(H_{w_4}  H_{w_5} H_{w_6})}.
\end{equation}

The theorems stated in the introduction are the translation to the analytic context of theorems \ref{teo-formula}, \ref{dobles} and \ref{simples}. In order to proof them,  we have only to check
that the bitangent lines associated to the odd two-torsion points $w_k\in J_\C[2]^{\rm odd}$ appearing in the theorems satisfy the hypothesis of the mentioned theorems. But this is immediate, using the analytic description (\ref{Weil}) of the Weil pairing.

\section{Some geometry around  the Frobenius formula}\label{frobenius}

While the closed solution to the Torelli problem given by theorem \ref{torelli} is satisfactory from the geometrical viewpoint,  it would be desirable  to have also a formula involving only Thetanullwerte instead of Jacobian Nullwerte.
We can derive this formula from our theorem just taking into account Frobenius formula, which we recall briefly.

Both the Thetanullwerte and the Jacobian Nullwerte can be interpreted as functions defined on the Siegel upper half space $\HH_3$: the idea consists in fixing the necessary characteristics and letting $Z$ run through $\HH_3$.  Frobenius formula relates certain Jacobian Nullwerte with  products of Thetanullwerte.

Three characteristics $\e_1, \e_2, \e_3$ are called asyzygetic if $e_2(\e_1-\e_2,\e_1-\e_3)=-1$. A sequence of characteristics is called asyzygetic if every triplet contained in it is asyzygetic.
A fundamental system  of characteristics is an asyzygetic  sequence $S=\{\e _{1},...,\e _{8}\}$ of  characteristics, with $\e_1,\e_2,\e_3$ odd and $\e_4,\dots,\e_8$ even.

\begin{theorem}[Frobenius-Igusa, \cite{Frobenius}, \cite{Igusa-odd}]
Let
$\e_1, \e_2, \e_3 $ be three odd characteristics, and take $w_i:=w_i(Z):=\e_i'+Z.\e_i''$. There is an equality    on $\HH_3$ of the form
 $$
 [w_1, w_2,
w_3](Z)=\pm\pi^3 \q[w_4](Z)\cdots\q[w_8](Z)
$$
  if and only if
$\e_1, \e_2, \e_3$ are asyzygetic, and in this case the characteristics $\e_4,\dots,\e_8$ corresponding to $w_4,\dots, w_8$ are the unique completion of $\e_1,\e_2, \e_3$ to a fundamental system.
\end{theorem}

The sign in the equality above can be determined for every particular fundamental system.

As we have mentioned after theorem \ref{teo-formula}, all the triplets of bitangent lines involved in the expression above are asyzygetic. This implies them that we can express the corresponding Jacobian Nullwerte as products of Thetanullwerte. We have the following identities:
$$
\begin{array}{ll}
\displaystyle
[w_1,w_2,w_3]=\pi^3 \prod_{k=1}^5\q[a_k](Z),\qquad
& \displaystyle[w_1',w_2',w_3']=-\pi^3 \prod_{k=1}^5\q[b_k](Z),\\
\displaystyle
[w_7,w_2,w_3]=\pi^3\prod_{k=1}^5\q[c_k](Z),\qquad &
\displaystyle[w_7,w_2',w_3']=\pi^3 \prod_{k=1}^5\q[d_k](Z),\\
\displaystyle
[w_1,w_7,w_3]=-\pi^3\prod_{k=1}^5\q[e_k](Z),\qquad &
\displaystyle[w_1',w_7,w_3']=-\pi^3 \prod_{k=1}^5\q[f_k](Z),\\
\displaystyle
[w_1,w_2,w_7]=\pi^3\prod_{k=1}^5\q[g_k](Z),\qquad &
\displaystyle[w_1',w_2',w_7]=\pi^3 \prod_{k=1}^5\q[h_k](Z),
\end{array}
$$
where
$$ \begin{array}{lll}
{}^ta_1={}^t(0,0,0)+{}^t(0,0,0).Z,\quad &
{}^ta_2={}^t(0,0,0)+{}^t(\frac12,0,0).Z,\quad &
{}^ta_3={}^t(\frac12,0,\frac12)+{}^t(\frac12,\frac12,\frac12).Z,\\
{}^ta_4={}^t(\frac12,\frac12,0)+{}^t(\frac12,\frac12,0).Z,\quad &
{}^ta_5={}^t(\frac12,\frac12,\frac12)+{}^t(\frac12,0,\frac12).Z,\\\\
{}^tb_1={}^t(0,0,0)+{}^t(\frac12,0,0).Z,\quad &
{}^tb_2={}^t(\frac12,0,0)+{}^t(0,0,0).Z,\quad &
{}^tb_3={}^t(\frac12,0,\frac12)+{}^t(\frac12,\frac12,\frac12).Z,\\
{}^tb_4={}^t(\frac12,\frac12,0)+{}^t(\frac12,\frac12,0).Z,\quad &
{}^tb_5={}^t(\frac12,\frac12,\frac12)+{}^t(\frac12,0,\frac12).Z,\\\\
{}^tc_1={}^t(0,0,0)+{}^t(0,0,0).Z,\quad &
{}^tc_2={}^t(0,\frac12,0)+{}^t(0,0,0).Z,\quad &
{}^tc_3={}^t(0,\frac12,\frac12)+{}^t(0,\frac12,\frac12).Z,\\
{}^tc_4={}^t(\frac12,0,0)+{}^t(0,\frac12,0).Z,\quad &
{}^tc_5={}^t(\frac12,0,\frac12)+{}^t(\frac12,\frac12,\frac12).Z,\\\\
{}^td_1={}^t(0,0,0)+{}^t(0,\frac12,0).Z,\quad &
{}^td_2={}^t(\frac12,0,0)+{}^t(0,0,0).Z,\quad &
{}^td_3={}^t(\frac12,0,\frac12)+{}^t(\frac12,\frac12,\frac12).Z,\\
{}^td_4={}^t(\frac12,\frac12,0)+{}^t(0,0,0).Z,\quad &
{}^td_5={}^t(\frac12,\frac12,\frac12)+{}^t(0,\frac12,\frac12).Z,\\\\
{}^te_1={}^t(0,0,0)+{}^t(0,0,0).Z,\quad &
{}^te_2={}^t(0,0,0)+{}^t(0,0,\frac12).Z,\quad &
{}^te_3={}^t(0,0,\frac12)+{}^t(0,\frac12,0).Z,\\
{}^te_4={}^t(\frac12,0,\frac12)+{}^t(\frac12,\frac12,\frac12).Z,\quad &
{}^te_5={}^t(\frac12,\frac12,\frac12)+{}^t(0,0,0).Z,\\\\
{}^tf_1={}^t(0,\frac12,\frac12)+{}^t(0,0,0).Z,\quad &
{}^tf_2={}^t(\frac12,0,0)+{}^t(0,0,0).Z,\quad &
{}^tf_3={}^t(\frac12,0,0)+{}^t(0,0,\frac12).Z,\\
{}^tf_4={}^t(\frac12,0,\frac12)+{}^t(0,\frac12,0).Z,\quad &
{}^tf_5={}^t(\frac12,0,\frac12)+{}^t(\frac12,\frac12,\frac12).Z,\\\\
{}^tg_1={}^t(0,0,0)+{}^t(0,0,0).Z,\quad &
{}^tg_2={}^t(0,0,0)+{}^t(0,\frac12,\frac12).Z,\quad &
{}^tg_3={}^t(0,0,\frac12)+{}^t(0,0,0).Z,\\
{}^tg_4={}^t(\frac12,0,\frac12)+{}^t(\frac12,\frac12,\frac12).Z,\quad &
{}^tg_5={}^t(\frac12,\frac12,0)+{}^t(0,0,\frac12).Z,\\\\
{}^th_1={}^t(0,\frac12,0)+{}^t(0,0,\frac12).Z,\quad &
{}^th_2={}^t(\frac12,0,0)+{}^t(0,0,0).Z,\quad &
{}^th_3={}^t(\frac12,0,0)+{}^t(0,\frac12,\frac12).Z,\\
{}^th_4={}^t(\frac12,0,\frac12)+{}^t(0,0,0).Z,\quad &
{}^th_5={}^t(\frac12,0,\frac12)+{}^t(\frac12,\frac12,\frac12).Z.
\end{array}
$$
Observe that there are some coincidences between the fundamental systems giving these identities: for instance,
$c_1=a_1$ and $c_5=a_3$. In fact, every two share two of their even characteristics. Thus, when we plug the equalities above into theorem \ref{torelli}, we can make some simplifications. We obtain

\begin{corollary}\label{torelli-theta} Let
$$
\begin{array}{l}
A_1=\theta[c_2](Z)\theta[c_3](Z)\theta[c_4](Z)\theta[d_1](Z)\theta[d_4](Z)\theta[d_5](Z),\\
A_2=\theta[e_2](Z)\theta[e_3](Z)\theta[e_5](Z)\theta[f_1](Z)\theta[f_3](Z)\theta[f_4](Z),\\
A_3=\theta[g_2](Z)\theta[g_3](Z)\theta[g_5](Z)\theta[h_1](Z)\theta[h_3](Z)\theta[h_4](Z).\\
\end{array}
$$
An equation  for $\C$, defined over $K$ up to a normalization, is:
$$
\C: (A_1 X_1 Y_1+A_2 X_2 Y_2-A_3 X_3 Y_3)^2-4 A_1 A_2 X_1 X_2Y_1 Y_2=0.
$$
\end{corollary}

\noindent
{\textbf{Remarks}:} We could have written a presentation of $\C$ using quotients of Thetanullwerte, so that the coefficients could be interpreted as modular functions for certain level congruence subgroups, as is the case for the coefficients of the formula in theorem \ref{torelli}. The advantage of the formula we have written is that it involves only 18 values of the Riemann theta function. Hence, from the computational point of view, this is a more satisfactory formula, since  the convergence of $\q$ is faster than the convergence of its derivatives, and moreover we need less evaluations.


We finish with a geometrical description of the fundamental systems  appearing in the Frobenius formula.
We have considered the general hyperelliptic case in \cite{Guardia-fourier}, and we now explain the situation for
non-hyperelliptic genus three curves:


\begin{theorem}
Let $X_1=0, X_2=0, X_3=0$ be three asyzygetic bitangent lines. Write
 $$
 S_{X_1 X_2}\cap S_{X_1,X_3}=\{ X_1=0, X_4=0, \dots, X_8=0\},
 $$
according to proposition \ref{steiner-asyzygetics}).
 Let
$$
\begin{array}{l}
W_j=\frac12\divi(X_j),\quad j=1,2,3;\\

W_j=\frac12\divi \left(\frac{X_2 X_3}{X_j}\right)=W_2+W_3-\frac12\divi(X_j),\quad j=4,\dots, 8.

\end{array}
$$
The odd two-torsion points $w_1,\dots, w_8\in J_\C$ corresponding to these  semicanonical divisors
through the Abel-Jacobi map (\ref{Abel-Jacobi})
 form a fundamental system.
\end{theorem}


\noindent {\bf Proof:} We must show that $e_2(w_1-w_i,w_1-w_j)=1$ for
all pairs $i,j\in\{1,\dots, 8\}$. By the construction of the $W_i$,
it is enough to see that $e_2(w_1-w_2,w_1-w_3)=1$,
$e_2(w_1-w_2,w_1-w_4)=1$ and $e_2(w_1-w_4,w_1-w_5)=1$. The first
equality is immediate, since $X_1=0, X_2=0, X_3=0$ are asyzygetic.
The remaining two are derived easily: we compute the Weil pairing applying the Riemann-Mumford relation \ref{Riemann-Mumford} with $D=W_1$. We have, for instance:
$$
\begin{array}{l}
\hskip-5mm e_2(w_1-w_2,w_1-w_4)=e_2(W_1-W_2,W_1-W_4)\\\\
=
h^0(W_1)+h^0(4W_1-W_2-W_4)+h^0(2W_1-W_2)+h^0(2W_1-W_4)\\\\
=1+h^0(2K-W_2-W_4)+h^0(W_2)+h^0(W_4)\\\\
=3+h^0(W_2+W_4)=1
\pmod2.
\end{array}
$$

\end{document}